\documentclass[11pt]{amsart}
\usepackage{amscd, amssymb, mathrsfs, 
dsfont, xcolor, 
amsmath, tikz-cd, soul}
\usepackage{kotex}
\usepackage[backref=page]{hyperref}   
\usepackage{collectbox}

\usepackage{enumitem}



\makeatletter
\@namedef{subjclassname@2020}{\textup{2020} Mathematics Subject Classification}
\makeatother

\input{xy}
\xyoption{all}

\newtheorem{Thm}{Theorem}[section]

\newtheorem{Prop}[Thm]{Proposition}
\theoremstyle{definition}

\newtheorem{Def/Thm}[Thm]{Definition/Theorem}
\newtheorem{Cor}[Thm]{Corollary}
\newtheorem{lemma}[Thm]{Lemma}

\newtheorem{Rmk}[Thm]{Remark}

\numberwithin{equation}{section}
\newcommand{\ti }{\times}
\newcommand{\ot }{\otimes}

\newcommand{\Hom }{{\mathrm{Hom}}}

\newcommand{\cA}{{\mathcal{A}}}
\newcommand{\cB}{\mathcal{B}}

\newcommand{\cN}{\mathcal{N}}

\newcommand{\cP}{{\mathcal{P}}}
\newcommand{\cQ}{{\mathcal{Q}}}

\newcommand{\GG }{{\mathbb G}}

\newcommand{\ZZ }{{\mathbb Z}}

\newcommand{\ka }{{\alpha}}

\newcommand{\Ch}{\mathrm{Ch}}

\newcommand{\End}{\mathrm{End}}

\newcommand{\Perf}{\mathrm{Perf}}
\newcommand{\Mod}{\mathrm{Mod}}
\newcommand{\Hoch}{C}

\newcommand{\MC}{\mathrm{MC}}

\newcommand{\id}{\mathrm{id}}

\begin{document}
\title{Chern characters for curved dg-algebras}  

\author[K.~Chung]{
Kuerak Chung}
\address{Department OF Mathematics Education\\
  Korea National University of Education\\
250 Taeseongtabyeon-ro, Gangnae-myeon\\
Heungdeok-gu, Cheongju-si, Chungbuk 28173\\
 Republic of Korea}
 \email{krchung@kias.re.kr}

\author[B.~Kim]{Bumsig Kim}
\address{Korea Institute for Advanced Study\\
85 Hoegiro, Dongdaemun-gu \\
Seoul 02455\\
Republic of Korea}
\email{bumsig@kias.re.kr}

\author[T.~Kim]{Taejung Kim}
\address{Department OF Mathematics Education\\
  Korea National University of Education\\
250 Taeseongtabyeon-ro, Gangnae-myeon\\
Heungdeok-gu, Cheongju-si, Chungbuk 28173\\
 Republic of Korea}
\email{tjkim@kias.re.kr}


\thanks{K.~Chung and T.~Kim are supported by NRF-2018R1D1A3B07043346. B.~Kim is supported by KIAS individual grant MG016404.}

\begin{abstract}
We construct a quasi-inverse of the cochain map on the negative cyclic complexes of the second kind induced from the quasi-Yoneda embedding on a curved dg algebra. 
This gives an explicit formula for the Chern character  of a perfect module. 
\end{abstract}

\subjclass[2020]{Primary 14A22; Secondary 16E40, 18G80}

\keywords{Matrix factorizations, Negative cyclic homology, Hochschild homology, Categorical Chern characters}
\maketitle


\section{Introduction}

\subsection{Overview}
Let $k$ be a field and let $\GG$ be either $\ZZ$ or $\ZZ/2$. 
Curved differential  $\GG$-graded $k$-algebras  (for short, cdg algebras) and their modules appear in relation with deformations of dg algebras \cite{KLN}, Koszul duals of nonhomogeneous quadratic algebras \cite{Pos}, mirror duals of Fano varieties in mirror symmetry \cite{Giv}, and LG-models in LG/CY correspondence \cite{Wit}, etc.

 Given a cdg algebra $\cA:=(A,d,h)$, its ordinary negative cyclic homology $HN_*(\cA)$ (and its Hochschild homology $HH_*(\cA))$ vanishes if its curvature $h$ is nonozero  (see \cite{CT, PP}). Thus, its negative cyclic homology of the second kind $HN^{II}_*(\cA)$ is used as a pseudo-equivalence (called ``quasi-Morita'') invariant of cdg algebras. Moreover, if $\cA $ is a cofibrant dg algebra $(A,d,h=0)$, then the natural comparison map $HN_*(\cA) \to HN_*^{II}(\cA)$ is an isomorphism \cite{PP}. 
 
On the other hand, there is another ``quasi-Morita" invariant of  $\cA$:  the Grothendieck group  $K_0(\cA)$  of 
the idempotent completion $[\overline{\Perf(\cA)}]$ for  the  homotopy category  $[ \Perf(\cA) ]$ of the dg category $\Perf (\cA)$ of perfect right $\cA$-modules.

Two  invariants are naturally related by the \emph{Chern Character map with values in the negative cyclic homology of the second kind}
\[ \Ch_{HN}^{II}:K_0(\cA) \to HN^{II}_0(\cA). \]
In this paper, we obtain a formula for $\Ch_{HN}^{II}(\cP) \in HN_0^{II}(\cA)$ which manifests the various known results
\cite{BW,Se,Shk: HRR} as special cases and is extended to a global geometric version for smooth curved algebras via Chern-Weil theory in \cite{CKK}.

 \subsection{Main results}

  Let $\cA=(A,d,h)$ be a cdg algebra. A right $\cA$-module $\cP=(P,\delta_P)$ is perfect if and only if $\cP$ is a direct summand of a twist $\cN_\alpha:=(N,d_F+\alpha)$ of shifted free right $\cA$-module $N=\oplus_{i=1}^{l}A[n_i]$ with $\alpha \in \End^1_A(N)$, $d_F(\alpha)+\alpha^2=\lambda_{-h}$. Here $\lambda_{-h}$ is the left multiplication by $-h$. Then, there are a non-unital dg algebra homomorphism $F:=i \circ (\;) \circ j:\End_\cA(\cP) \to \End_\cA(\cN_\alpha)$, a cdg algebra isomorphism $(\id_P,\alpha):\End_\cA(\cN_\alpha) \to \End_\cA(\cN_0)$ inducing a cochain map $(\id_P,\alpha)_* \circ F_*(u):$
  \[\overline{HN}^{II}(\End_\cA(\cP)) \to \overline{HN}^{II}(\End_\cA(\cN_\alpha)) \to \overline{HN}^{II}(\End_\cA(\cN_0)) \] 
and the Chern character of $\cP$ is given by
\[ \Ch^{II}_{HN}(\cP) = \mathrm{Tr}((\id_P,\alpha)_*( F_*(u) ([1_P]))) \in \overline{HN}^{II}_0(\cA)\]
where $\mathrm{Tr}:\overline{HN}^{II}(\End_\cA(\cN_0)) \to \overline{HN}^{II}(\cA)$ is the generalized trace map; see \S~\ref{sub: gen tr} for the definition.
  
More concretely, in \S~\ref{sub: via eta} we show that for a perfect $\cA$-module $\cP$ given by an idempotent $\pi$ on $\cN_\alpha$, $\Ch ^{II}_{HN}
 (\cP)$ is represented by a cocycle in the normalized negative cyclic complex of the second kind $(\overline{C}^{II}(\cA)[[u]],\bar{b}+u\overline{B})$:

\begin{equation} 
 \sum_{j=0}^{\infty } (-1)^j  \mathrm{Tr}(\pi[\alpha^j])+\sum_{n \ge 1}\sum_{ (j_0, ..., j_{2n}) \in \ZZ ^{2n+1}_{ \ge 0}} c_{n,J} \mathrm{Tr}  ((2\pi-\  1_N ) [\underbrace{\alpha^{j_0}|\pi|\alpha^{j_1}|\pi|\cdots|\pi|\alpha^{j_{2n}}}_{2n+ J} ] )u^n 
 \end{equation}
where $J= \sum_{k=0}^{2n} j_k $, $c_{n, J} 
:=(-1)^{n+J } \frac{(2n)!}{2(n!)}$, $\ka ^ r = \underbrace{  \ka | ... | \ka }_{r}$ for $r= j , j_0, ..., j_{2n}$. Moreover, if $h=0$ and $\GG=\ZZ$, then the  Chern character
$\Ch  _{HN} (\cP)$ of the first kind is represented by the same formula above but with a finite sum. In \S~\ref{sub: homotopy}, we briefly mention a Chern character formula for a homotopy direct summand of $\cN_\alpha$. 

\subsection{Acknowledgements} 
 K.C. and T.K. thank Chanyoung Sung and Hoil Kim for encouragement and support. The second author, Professor Bumsig Kim, passed away in the final stages of preparing this manuscript. We would like to extend our deepest gratitude to him for showing his professional and personal generosity and kindness.

 \subsection{Conventions and notation}

This paper is a sequel of \cite{CKK}. Unless otherwise stated, 
we freely use the conventions from \cite{CKK} and  the definitions and notation from Section 2 of \cite{CKK}, in particular for the category of mixed complexes, 
the (normalized) Hochschild chain complexes (of the second kind), the (normalized) Hochschild homology (of the second kind), 
the (normalized) mixed Hochschild chain complexes (of the second kind), the (normalized) negative cyclic chain complexes (of the second kind), 
the (normalized) negative cyclic homology (of the second kind) etc. We refer to \cite{She} for similar definitions and properties for $A_\infty$-categories.

\section{Categorical Chern characters}\label{sec: cat chern}

 \subsection{Definitions of Chern characters} 
 
 \subsubsection{} Given a small cdg category $\cA$, 
 let $q\Perf(\cA)$ be the closure of $\cA$ in the cdg category $q\Mod(\cA)$ of right quasi-modules over $\cA$ under shift, direct sum, twist, and passage to a direct summand.
 The objects of $q\Perf(\cA)$ are called {\em perfect right quasi-modules} over $\cA$. 
Let $\Perf(\cA)$ be the full subcategory of $q\Perf (\cA)$ consisting of all right $\cA$-modules, called \emph{perfect right $\cA$-modules}. 
 
 \subsubsection{}  The homotopy category $[\Perf(\cA)]$ of $\Perf(\cA)$ is naturally triangulated.
Let $ K_0([\Perf(\cA)]) $ be the Grothendieck group of the triangulated category $[\Perf(\cA)]$.
By the tautological assignments $P \mapsto [1_P]$, we obtain a homomorphism $taut ^{II} : K_0 ([\Perf(\cA)]) \to  \overline{HN}^{II}_0(\Perf(\cA)) $.
Together with induced maps from natural embeddings, $taut ^{II}$ makes the following diagram define Chern character maps (see \cite{BW}):
\[\xymatrix{ K_0([\Perf(\cA)]) \ar[r]^(.49){taut^{II}}  \ar[rd]_{taut} \ar@/^2pc/[rrr]^{ \Ch_{HN}^{II}} \ar@/_4pc/@{-->}[rrd]_{ \Ch_{HN}}   
                           & \overline{HN}^{II}_0(\Perf(\cA)) \ar[r]^{\cong}  & \overline{HN}^{II}_0(q\Perf(\cA))  & \ar[l]_(.4){\cong}  \overline{HN}^{II}_0(\cA) \\
                                                                 & \overline{HN}_0(\Perf (\cA))  \ar[u] &  \ar@{-->}[l]_(.4){\cong}  \overline{HN}_0(\cA)  &  } \]
where the dotted arrows make sense when $\cA$ is a dg category.

\subsection{Hochschild morphism for a semifunctor}
For an arbitrary  non-unital functor $F : \cA \to \cB$ between dg categories, there is a functorial cochain map $F_* :  (C (\cA ), b)  \to (C (\cB) , b)$.
We will construct its extension 
\[ F_* (u)  : (C (\cA ) [[u]], b + uB)  \to (C (\cB) [[u]], b + uB) \]
by letting $F_*(u) := p(u) \circ F_*^e \circ \iota (u)$; see \S~\ref{sub: non unit} for the definition of $F_*^e$ and \S~\ref{subsub: iu pu} for the definitions of  $\iota (u), p(u)$.

\subsubsection{}\label{sub: non unit}  In \S~\ref{sub: non unit} we allow that  $\mathcal{A}$ is a dg category possibly  without identities. 
 Denote the dg category obtained from $\mathcal{A}$ by adjoining identities by $\mathcal{A}^+$; its objects are the same objects as $\mathcal{A}$ and its morphisms are defined by
\[ \mathcal{A}^+(x, y)=\mathcal{A}(x, y) \;\; \mathrm{if} \;\; x \ne y, \quad \mathcal{A}^+(x, x)=\mathcal{A}(x, x)\oplus k\cdot e_x, |e_x|=0. \]
The composition on $\mathcal{A}^+$ is defined in the obvious way. Consider the subcomplex $D$ of $(C(\mathcal{A}^+),b^+)$ generated by 
$e_x$ and 
$a_0[a_1| \cdots|a_n]$ with $a_i=e_{x}$ for some $i>0$ and any $x \in \cA$.
We denote $(C(\mathcal{A}^+)/D,\overline{b^+})$  by $(C^e(\mathcal{A}),b^e)$ and  call it  the \emph{non-unital Hochschild complex} 
of $\cA$. Let  $\mathrm{cod}(a_i)$ denote  the codomain of $a_i$ and define Connes' differential $B^e:C^e(\mathcal{A}) \to C^e(\mathcal{A})$ by 
\[B^e(a_0[a_1|\cdots |a_n]):=e_{\mathrm{cod}(a_0)}[a_0|a_1|\cdots |a_n] \] 
so that $\mathrm{MC}^e(\cA):=(C^e(\mathcal{A}),b^e,B^e)$ is a mixed complex and  called the the \emph{non-unital Hochschild mixed complex} of $\cA$.

\begin{Rmk} The following ones hold; see  \label{rmk: semi hom}  \cite[ \S~3.5]{She}.
\begin{enumerate}[label=(\roman*)]
\item
The non-unital Hochschild complex of $\mathcal{A}$ is, by definition, the reduced Hochschild complex of $\mathcal{A}^+$; see \cite[\S~1.4.2]{Loday}.
\item\label{item semi Hoch}
For a dg category $\cA$ with identities, the composition of the natural maps
\begin{equation}\label{def: j}
j:(C(\mathcal{A}),b) \hookrightarrow (C(\mathcal{A}^+),b^+) \to (C^e(\mathcal{A}),b^e) 
\end{equation}
is a quasi-isomorphism.

\item\label{item: semi hom}
Any semifunctor $F:\mathcal{A} \to \mathcal{B}$  induces a canonical dg functor
\[F^+:\mathcal{A}^+ \to \mathcal{B}^+ \]
and a morphism of mixed complexes
\[ F^e:  (C^e(\cA) , b^e,B^e) \to (C^e(\mathcal{B}),b^e,B^e)\]
compatible with the quasi-isomorphisms in item~\ref{item semi Hoch}.
\item
The construction of non-unital Hochschild complex is functorial; for semifunctors $F;\mathcal{A} \to \mathcal{B}, G:\mathcal{B} \to \mathcal{C}$,
\[ (G \circ F)^e=(G^e\circ F^e).\]
\end{enumerate}
\end{Rmk}

\subsubsection{}\label{subsub: iu pu}
Let $\mathcal{A}$ be a dg category with identities as usual. 
The map $j: (C(\mathcal{A}),b) \to (C^e(\mathcal{A}),b^e)$ in \eqref{def: j} is not compatible with Connes' operators,
 but we can still construct homotopy inverses between $(C(\cA)[[u]],b+uB)$ and $(C^e(\cA)[[u]],b^e+uB^e)$ as $k[[u]]$-modules which  extends $j$. 
 This construction was done by Shklyarov (see \cite{Shk: NC}) as follows: 
 
As graded $k$-modules
\[ C^e(\mathcal{A})=C(\mathcal{A})\oplus C^+(\mathcal{A}) \]
where $C^+(\mathcal{A})=\{a' \in C^e(\mathcal{A})| a'=e_{\mathrm{cod}(a_1 ) }[a_1|\cdots|a_n], n \ge 1 \}$ and  $j(a)=(a,0)$. 
Consider the following $k$-linear maps:
\begin{align*}  \iota(u):(C(\mathcal{A})[[u]],b+uB) & \to (C^e(\mathcal{A})[[u]],b^e+uB^e) \\ 
                a &  \mapsto a+us^esN(a)  \end{align*}
                and 
\begin{align*}   p(u):(C^e(\mathcal{A})[[u]],b^e+uB^e) & \to (C(\mathcal{A})[[u]],b+uB) \\
       (a,a') & \mapsto a+(1-t^{-1})\mu(a') ,  \end{align*}
where $\mu:C^+(\mathcal{A}) \to C(\mathcal{A}) ;  \  e_{\mathrm{cod}(a_1)} [a_1|\cdots |a_n] \mapsto 1_{\mathrm{cod}(a_1) } [a_1|\cdots |a_n]$.

For a dg algebra $A$, $\iota(u), p(u)$ are shown to be homotopy inverses over $k((u))$ in \cite[ \S~3.2]{Shk: NC} and the same proof for a dg category $\mathcal{A}$ works over $k[[u]]$ as follows:
\begin{Rmk}\label{rmk: iu pu}  
  \hfill
  \begin{enumerate}[label=(\roman*)]
\item
$\iota(u), p(u)$ are $k[[u]]$-linear morphisms.
\item
$\iota(u), p(u)$ are homotopy inverses as $k[[u]]$-morphisms:
$p(u)\circ \iota(u)$ is homotopic to $\mathrm{id}_{C(\mathcal{A})[[u]]}$ with a homotopy $H(u):=u(1-t^{-1})sssN$ and $\iota(u)\circ p(u)$ is homotopic to $\mathrm{id}_{C^e(\mathcal{A})[[u]]}$ with a homotopy $H^e(u):(a,a') \to (0,s^e\mu(a'))$.
\end{enumerate}
\end{Rmk}

\subsection{A lemma}

Suppose that $P$ is a direct summand of $N$ in a dg category $\cA$. In other words, there are degree $0$ closed homomorphisms $i : P \to N$ 
and $j: N \to P$ such that $1_P = j \circ i$ and denote $\pi : = i\circ j $. Then, we have a non-unital functor $F : \{ P, N \} \to \{ N \}$ of full dg subcatgories of $\cA$ defined by
\begin{equation}\label{eqn: s-fun} f \mapsto    i\circ f \circ j, \quad  f ' \mapsto i \circ f', \quad f'' \mapsto f '' \circ j, \quad g \mapsto g \end{equation}
 for $f \in \End P, \; f' \in \Hom (N, P), \; f'' \in \Hom (P, N), \; g \in \End (N)$.

In $(\Hoch \{ P \}  [[u]] , b + uB)$, 
let   \begin{align*} \gamma_P & :=   1_P +   \sum _{i=1}^{\infty}        (-1)^{i}  \frac{(2i)!}{2(i!)}  2 \cdot 1_P  [\underbrace{1_P |  \cdots | 1_P }_{2i} ] u^{i} ; \text{ and }  \\
  \eta _{\pi} & :=      \pi +   \sum _{i=1}^{\infty}        (-1)^{i}  \frac{(2i)!}{2(i!)}  (2\pi -1_N)  [\underbrace{\pi|  \cdots | \pi}_{2i} ] u^{i}  \\
   & =  \pi - (2\pi - 1_N) [ \pi | \pi] u + 6(2 \pi - 1_N)[ \pi| \pi | \pi | \pi]  u^2 + \cdots . \\
 \end{align*}

\begin{lemma}\label{lem: general eta}
There is a cochain map  $F_*(u) :  \Hoch  \{ P, N \}  [[u]] \to \Hoch  \{  N \}  [[u]] $ between the negative cyclic complexes
 such that:
 \begin{enumerate} 
 \item\label{item inv}  The composition $F_*(u) \circ inc_* $ is the identity in the homology level $HN_* \{ N \}   \to HN_* \{ N \}$.
\item\label{item u0}  When $u=0$, it is $F_*$, i.e., $F_* (u) |_{u=0} = F_*$. 
\item\label{item eta}  Let $quot: \Hoch  \{  N \}  [[u]] \to \overline{C} \{  N \}  [[u]]$ denote the quotient map. Then 
$quot \circ F_*(u) (\gamma_P) = \eta _{\pi} $.
\end{enumerate}
In particular $[1_P]= [\eta _{\pi}]$ in $\overline{HN}_* \{ P, N \}$.
\end{lemma}
\begin{proof}
Since the semifunctor $F$ above does not preserve the unit $1_P$, we need a new idea. 
There are a so-called non-unital mixed Hochschild  complex $(C^{e} (\cA), b^e, B^e)$
and cochain maps 
\begin{align*} \iota (u) : &  (C (\cA )[[u]] , b+ uB) \to (C^{e} (\cA) [[u]] , b^e + uB^e) ; \text{ and }\\
  p(u): & (C^e (\cA) [[u]] , b^e + uB^e) \to (C (\cA) [[u]], b + uB), \end{align*}
which are homotopy inverses to each other; see Remark~\ref{rmk: iu pu}.
The semifunctor $F$ induces a cochain map 
\[ F_*^e : (C^e\{ P, N\}[[u]], b^e + u B^e)  \to (C^e \{N\}[[u]], b^e + uB^e) \] as explained in 
Remark~\ref{rmk: semi hom}~\ref{item: semi hom} such that $F_*^{e} |_{C^e \{ N \} [[u]] } = \mathrm{id}_{C^e \{ N \} [[u]] } $. 
Let $F_*(u) := p(u) \circ F_*^e \circ \iota (u)$. Then \eqref{item inv} and \eqref{item u0} are clear. 
Item \eqref{item eta} is a straightforward computation from the definitions of $\iota (u)$ and $p(u)$. 
\end{proof}

\begin{Cor}\label{prop: eta}  Two cycles $1_P$ and $\eta _{\pi}$ are  homologous in $(\overline{C} \{ P, N \} [[u]], b+ uB)$.
\end{Cor} 
We note that Corollary~\ref{prop: eta} is proved in \cite{CKK} by a very different method.

\subsection{A Chern character formula}\label{sub: via eta}

Let $\cA=(A, d, h)$ be a cdg algebra and let $\cP =(P,\delta_P) \in \Perf (\cA)$. 
Recall that there is a canonical isomorphism
\begin{equation}\label{eqn: can iso}  HN_*^{II} (\Perf (\cA)) \cong HN _*^{II} (\cA)  \end{equation}
induced by the embeddings $\cA \to q\Perf (\cA)$ and $\Perf (\cA) \to q\Perf (\cA )$.
In this subsection we apply Corollary~\ref{prop: eta}
to get a formula for $\Ch^{II}_{HN} (\cP)$, which works as well as for   $\Ch_{HN} (\cP)$ if $h=0$ and $\GG = \ZZ$.

Let $\cN_0: = (N, d_F)$ denote a finitely generated free cdg module over $\cA=(A, d, h)$, i.e.,  a finite sum of $A$'s up to degree shifts with 
the induced differential $d_F$ from the differential $d$ of $A$. We call $\cN_0$ a \emph{finitely generated free quasi-module} over $(A, d,h)$. 
Let $\cP=(P,\delta_P)$ be a direct summand of twisted $\cA$-module $\cN_{\ka}:=(N, d_F + \ka)$ of $\cN_0$ with $\ka \in \End ^1_A (N)$. Note that $d_F (\ka) + \ka ^2 = \lambda _{-h}$ where $\lambda _{-h}$ is the left multiplication by $-h$. The module $\cN_{\ka}$ is called a
 \emph{finitely generated semi-free module} over $\cA$ and note that any perfect $\cA$-module $\cP$ is represented as a direct summand of a finitely generated semi-free $\cA$-module. Consider the cdg isomorphism  \[ (\mathrm{id}_N, \ka) : \{ \cN_{\ka} \} \to \{ \cN_0 \}  \] 
between full cdg subcategories of $q\Perf (\cA)$ consisting of the indicated object only. For example,  $\{\cN_{\ka}\}$ is the dg algebra $\End_\cA(\cN_\ka):=(\End_A(N),d_{\cN_\alpha}=[d_F+\alpha,\;])$.
Then we get a Chern character formula for $\cP \in \Perf (\cA)$ as in the following proposition.

\begin{Prop}\label{prop: via eta} Let $\cP \in \Perf (A,d,h)$ be a direct summand of $\cN_{\ka}=(N,d_F+\alpha)$ given by a closed idempotent $\pi:N \to N$.
\begin{enumerate}
\item Under \eqref{eqn: can iso}, $ \Ch ^{II}_{HN} (\cP) $ is representable by a cocycle 
\[ \mathrm{Tr}  ((\mathrm{id}_{N}, \ka )_* (\eta _{\pi} ))  \in \overline{C}^{II} (A, d, h)[[u]] .\]
Here $\mathrm{Tr}$ is the $k[[u]]$-linear extension of the generalized trace map 
\begin{equation}\label{eqn: gen Tr} \overline{C}^{II} \{ \cN_0 \}\to \overline{C}^{II}  (A, d, h)  \end{equation} defined in \S~\ref{sub: gen tr}.

\item Suppose that $\GG = \ZZ$ and $h=0$. We may write $N=\oplus _{i=1}^l A[n_i]$ for some integers $n_i$, $l \ge 0$ for which
$\alpha : N \to N$ is a strictly upper triangular $l \times l$-matrix of entries
\[ \alpha_{ij} \in \mathrm{Hom}^1_{A}(A[n_j],A[n_i]) . \]
Then under the canonical isomorphism $HN_*(\Perf (A, d)) \cong HN _*(A, d) $, $ \Ch_{HN} (\cP) $ is representable by a cocycle 
\[ \mathrm{Tr}  ((\mathrm{id}_{N}, \ka )_* (\eta _{\pi} ))  \in \overline{C} (A, d)[[u]] .\]
\end{enumerate}

\end{Prop}

\begin{proof}  In this proof we write simply $\cA$ for the right quasi-module $(A,d)$ over $(A, d, h)$. (1) The cdg functor $(\mathrm{id}_N, \ka)$ can be extendable to a cdg functor 
$\widetilde{(\mathrm{id}_N, \ka)}: \{ \cN_{\ka}, \cN_0, \cA \} \to \{ \cN_0 , \cA \}$ which is a left inverse 
of the inclusion $\{ \cN_0 , \cA\} \to  \{ \cN_{\ka}, \cN_0, \cA \} $. Hence for every $x \in C^{II}\{ \cN_{\ka} \}[[u]]$ we note that 
$ x  $  and $ (\mathrm{id}_N, \ka )_* (x)$ are homologous in  $\overline{C}^{II} \{ \cN_{\ka}, \cN_0, \cA \}[[u]]$.
The map $\mathrm{Tr}$ in \eqref{eqn: gen Tr} can be extendable to a morphism 
$\widetilde{\mathrm{Tr}}: \overline{\MC}^{II} \{ \cN_0,  \cA \}\to \overline{\MC}^{II}  (\cA)  $ which is a left inverse
of the inclusion morphism of mixed complexes; see \S~\ref{sub: gen tr} for the extension. Again note that $z$ and $\mathrm{Tr} (z)$ are homologous 
in $\overline{C}^{II}  \{ \cN_0,  \cA \}  [[u]] $ for any $z \in \overline{C}^{II} \{ \cN_0 \}  [[u]]$.
Now by Corollary~\ref{prop: eta} and letting $x = \eta _{\pi}$ and $z= (\mathrm{id}_N , \ka )_* (\eta _{\pi})$ we conclude the proof.

(2) We have shown that  $\widetilde{\mathrm{Tr}} \circ \widetilde{(\mathrm{id}_N, \ka)}_* $ is a left inverse of the inclusion
 $\overline{\MC}^{II}(\cA) \to \overline{\MC}^{II} \{ \cN_{\ka}, \cN_0, \cA\}$. Since $\ka$ is strictly upper triangular, we can apply \eqref{tr 0}. Therefore
 $\widetilde{\mathrm{Tr}} \circ \widetilde{(\mathrm{id}_N, \ka)}_* $
  restricted to $\overline{\MC} \{ \cN_{\ka}, \cN_0, \cA \}$ lands in  $\overline{\MC} (\cA) $, which establishes the proof.
\end{proof}

\subsubsection{}\label{sub: gen tr}
Let $L$ be a graded $k$-module. For example $L$ is a finite sum of shifted $A$'s.
For positive integers $m_1$, $m_2$, $n$,  let $\mathrm{Mat}_{m_1\ti m_2} (L^{\ot n})$ denote the $k$-module of  all the $m_1 \ti m_2$ matrices with entries in  $L^{\ot n}$.
There is  the matrix multiplication 
\begin{align*} \bullet:   \mathrm{Mat}_{m_1\ti m_2} (L^{\ot _k n_1}) \ot _k  \mathrm{Mat}_{m_2\ti m_3} (L^{\ot_k n_2})  & \to \mathrm{Mat}_{m_1 \ti m_3} ( L ^{\ot _k { (n_1 + n_2) } }) 
 \\ (a_{ij}) \ot (b_{ij}) & \mapsto (\sum _j a_{ij} \ot b_{j k} ) \end{align*}  by the tensor algebra of $L$. 
If $ \mathrm{tr}:  \mathrm{Mat}_{m \ti m} ( L ^{\ot _k n} ) \to L ^{\ot _k n} $ denotes the supertrace map, then for $\phi _i \in \mathrm{Mat}_{m_i \ti m_{i+1}} (L)$
$i=0, ..., n$, $m_{n+1 }= m_0$, we define  the \emph{generalized
trace} map $\mathrm{Tr}$ by letting
\[ \mathrm{Tr} (\phi _0 \ot  \phi _1 \ot ... \ot \phi _n  ) := \mathrm{tr} (\phi _0 \bullet \phi _1 \bullet ... \bullet \phi _n ) . \]
 This is the graded version of the generalized trace map of  \cite[\S~1.2.1]{Loday}. Note that for $ m \ge m_{i}$
\begin{equation}\label{tr 0} \mathrm{Tr} ( \phi _0 \ot ...  \ot \phi _{i-1} \ot  \phi ^{\ot m}   \ot \phi _{i} \ot   ... \ot \phi _{n} ) =0 \end{equation} 
 if $\phi$ is a strictly upper triangular square matrix in $\mathrm{Mat} _{m_i \ti m_{i} }  (L)$.

\subsubsection{}
  Let $\cA=(A,d,h)$ be a cdg algebra.
Let $qF(\cA)$ denote the full subcategory of $q\Perf \cA$ consisting of finitely generated free quasi-modules
and let $sF (\cA)$ denote  the full subcategory of $\Perf \cA$ consisting of finitely generated semi-free modules.
The proofs of Lemma~\ref{lem: general eta} and Proposition~\ref{prop: via eta} show  there are  quasi-isomorphisms 

\begin{equation}\label{diag: perf A inv}
 \xymatrix{     &  \overline{C}^{II} \{ qF (\cA) , sF (\cA) \} [[u]] \ar@/_1.5pc/[dl]^(.62){\quad\  \{ (\mathrm{id} , 0)_*, (\mathrm{id}, \ka )_* \}  }  \ar@{^{(}->}[r] 
&   \overline{C}^{II} (q\Perf \cA )  [[u]]    \\
  \ar@<2pt>[d]^{\mathrm{Tr} \ }  \overline{C}^{II} (qF (\cA) ) [[u]] \ar@/^1.5pc/@{^{(}->}@<4pt>[ur] 
&  \ar@<2pt>[l]^{\{ (\mathrm{id}, \ka )_* \}  }   \overline{C}^{II} (sF (\cA) ) [[u]] \ar@{^{(}->}@<2pt>[r]   \ar@{_{(}->}[u]
&  \overline{C}^{II} (\Perf \cA ) [[u]] \ar@{_{(}->}[u]   \\
\overline{C}^{II} (\cA) [[u]] \ar@{^{(}->}@<2pt>[u]    &    C^{II} (sF (\cA) ) [[u]] \ar[u]^{quot}  \ar@{^{(}->}@<2pt>[r]  
&  \ar@<2pt>[l]^{ F_*(u)  }   C^{II} (\Perf \cA ) [[u]]     \ar[u]^{quot}   . } \end{equation}
All inclusions  are induced from the embeddings. The cochain maps $\{ (\mathrm{id}, \ka )_* \} $ and $\{  (\mathrm{id}, 0 )_*  , (\mathrm{id}, \ka )_*\} $ are induced from cdg functors. The right bottom map 
$F _* (u) $ is  induced from a semifunctor and non-unital mixed Hochschild complexes. 
All maps except the right bottom map are  induced from morphisms between mixed complexes. 
The cochain maps $\mathrm{Tr}$ and $\{ (\mathrm{id}, 0)_*, (\mathrm{id}, \ka)_* \} $ are left inverses of the corresponding inclusions, respectively.
Diagram~\eqref{diag: perf A inv} commutes in the homology level.  
Therefore the composition $\mathrm{Tr} \circ \{ (\mathrm{id}, 0)_*\} \circ quot \circ F_* (u) $ fits in a commutative diagram

\begin{equation}\label{diag: perf A inv hom}  \xymatrix{ \overline{HN}^{II}  (\cA) \ar[rr]^(.45){canonical} & &  \overline{HN}^{II} (q\Perf \cA ) & & \ar@/_2pc/[llll] HN ^{II} (\Perf \cA )  \ar[ll]_{canonical} . } 
\end{equation}
Thus 
  for $\cP$ in Proposition~\ref{prop: via eta}, $\Ch ^{II} _{HN} (\cP)$ is represented by a cocycle in $\overline{C}^{II}(\cA)[[u]]$:

  \begin{equation}\label{eqn: form II}   
  \sum_{j=0}^{\infty } (-1)^j  \mathrm{Tr}(\pi[\alpha^j])+\sum_{n \ge 1}\sum_{ (j_0, ..., j_{2n}) \in \ZZ ^{2n+1}_{ \ge 0}} c_{n,J} \mathrm{Tr}  ((2\pi-\  1_N ) [\underbrace{\alpha^{j_0}|\pi|\alpha^{j_1}|\pi|\cdots|\pi|\alpha^{j_{2n}}}_{2n+ J} ] )u^n 
 \end{equation}
where $J= \sum_{k=0}^{2n} j_k $, $c_{n, J} 
:=(-1)^{n+J } \frac{(2n)!}{2(n!)}$ and $\ka ^ r = \underbrace{  \ka | ... | \ka }_{r}$ for $r= j , j_0, ..., j_{2n}$. 
This proves formula (5.24) of \cite{BW}. 
When  $u$ is specialized to $0$, this  recovers  \cite[Theorem 2.14]{Se}. 

When $\GG = \ZZ$ and $h=0$,  the composition $\mathrm{Tr} \circ \{ (\mathrm{id}, \ka)_* \} \circ quot \circ F_* (u)$ restricted to  $C (\Perf \cA ) [[u]]$ 
is a quasi-inverse of the cochain map  $C(\cA)[[u]] \to C(\Perf \cA)[[u]]$ induced from the Yoneda embedding $\cA \hookrightarrow \Perf \cA$.   
Therefore diagram~\eqref{diag: perf A inv hom} without the superscripts $II$ makes sense and commutes. 
Hence
for $\cP$ in Proposition ~\ref{prop: via eta} (2),  $\Ch  _{HN} (\cP)$ is represented by a cocycle in $\overline{C}^{}(\cA)[[u]]$:

\begin{equation}\label{eqn: form I}  
   \sum_{j=0}^{l-1} (-1)^j \mathrm{Tr}(\pi[\alpha^j])+\sum_{ n\ge 1 }\sum_{0 \le j_0, ..., j_{2n}  \le l-1} c_{n,J} \mathrm{Tr}  ( (2 \pi-\ 1_N) [\underbrace{\alpha^{j_0}|\pi|\alpha^{j_1}|\pi|\cdots|\pi|\alpha^{j_{2n}}}_{2n+ J} ])u^n .
 \end{equation}
When  $u$ is specialized to $0$, this  recovers  \cite[Theorem 1.1]{Shk: HRR}.

\subsubsection{}\label{sub: homotopy}

  Let  $\cQ=(Q,\delta_Q) \in  \overline{\Perf(\cA)}$ be represented by a homotopy direct summand of a finitely generated semi-free $\cA$-module $N_\alpha=(N,d_N+\alpha)$. Then, there are closed $A$-module homomorphisms $i:Q \to N, \; j:N \to Q$  such that $j \circ i = \id_Q+[\delta_Q,H_1], H_1 \in \End_A^{-1}(Q), i \circ j = \pi+ [d_N+\alpha,H_2], H_2 \in \End_A^{-1}(N)$. 
 
 One can check that the following family $F:=\{F_n\}_{n \ge 1}:\{ \cN_\alpha, \cQ \} \to \{\cN_\alpha \}$ defines a non-unital  $A_\infty$-functor. For  $f_i \in \End_A(Q)$ with $\mathrm{cod}(f_1)=\mathrm{dom}(f_n)=\cQ$,
\[ F_n(f_1 \otimes \cdots \otimes f_n):=i\circ f_1 \circ H_{i_1} \circ f_2 \circ H_{i_2} \cdots \circ H_{i_{n-1}} \circ f_n \circ j \] 
where $H_{i_j}$ is $H_1$ or $H_2$ depending on $f_{i_j}$. For others, it is defined in the obvious way via \eqref{eqn: s-fun}. 
  This extends to a unital $A_\infty$-functor $F^e=\{F_n^e\}_{n \ge 1}:\{\cN_\alpha, \cQ\} \to \{\cN_\alpha\}$ and with this $F$, Lemma ~\ref{lem: general eta} is true (see \cite{She} for the definitions of $A_\infty$-category and $A_\infty$-functor, and for similar above facts about $A_\infty$-category). Thus,  
\[ \Ch^{II}_{HN}(\cQ)=\mathrm{Tr}\big((\id,\alpha)_*(F_*(u)(\gamma_Q))\big) \]
 where $F_*(u):=p(u)\circ F^e_* \circ i(u)$.

\end{document}